\documentclass[
12pt,
authoryear,
times
]{elsarticle}

\usepackage{amsmath}

\usepackage[colorlinks=true, %pdfstartview=FitV, linkcolor=black, citecolor= black, urlcolor= blue
]{hyperref}

\usepackage{lineno}

\journal{Acta Astronautica}

\begin{document}

\begin{frontmatter}

\title{ On Inclination Resonances in Artificial Satellite Theory\tnoteref{KePASSA} }
\tnotetext[KePASSA]{Preliminary results were presented in DyCoSS 2014, Rome, Italy, March 23--25, 2014}

\author%[roa]
{ Martin Lara\fnref{SDG} %\corref{cor1} 
}
\fntext[SDG]{ Space Dynamics Group, Polytechnic University of Madrid -- UPM, Madrid, Spain }

\ead{mlara0@gmail.com}

\address%[roa]
{Columnas de H\'ercules 1, ES-11000 San Fernando, Spain}


\date{}

\begin{abstract}
The frozen-perigee behavior of elliptic orbits at the critical inclination is usually displayed after an averaging procedure. However, this singularity in Artificial Satellite Theory manifests also in the presence of short-period effects. Indeed, a closed form expression relating orbital inclination and the ratio anomalistic draconitic frequencies is derived for the main problem, which demonstrates that the critical inclination results from commensurability between the periods with which the radial and polar variables evolve in the instantaneous plane of motion. This relation also shows that the critical inclination value is slightly modified by the degree of oblateness of the attracting body, as well as by the orbit's size and shape.
\end{abstract}


\end{frontmatter}

%\linenumbers

\section{Introduction}

The critical inclination problem of Artificial Satellite Theory (AST), a particular orbital inclination that ``freezes'' the perigee of elliptic orbits, has been qualified as ``the most celebrated problem in AST'' \citep{Jupp1988}. This special inclination of $63.4$
degrees (resp.~$116.6$ deg.) becomes apparent from the simple derivation of Lagrange planetary equations when the disturbing function of the main problem of AST, which is obtained after neglecting all harmonic coefficients in the geopotential except for the second order zonal harmonic coefficient ($J_2$), is averaged over the mean anomaly.
Thus, for instance, referring to the thorough derivations in Sec.~10.6 of Battin's celebrated book \cite{Battin1999}, it is simple to show  that the semi-major axis, eccentricity and inclination remain constant on average, and that there exists a constraint between these constant parameters and the mean frequencies of the motion. In particular, the mean orbital inclination can be explicitly written as a function of the ratio of the mean rate of rotation of the line of apsides to the orbital mean motion.
\par

It must be noted that the classical approach in Battin's book limits to the constant term of the expansion of the main problem disturbing function as a Fourier series in the mean anomaly, in this way missing important effects of the second order of $J_2$. Indeed, for a given mean semi-major axis, second order effects of $J_2$ reduce the number of elliptic frozen orbits at the critical inclination to just four isolated solutions with arguments of the perigee $0$, $\pi/2$, $\pi$, and $3\pi/2$, respectively, two of which are stable and the other two unstable \citep{Hagihara1961,Kozai1961}. Besides, other second order effects of the artificial satellite problem may introduce qualitative and quantitative changes in the phase portrait of orbits at the critical inclination, as, for instance, when the number of zonal harmonics included in the truncation of the gepotential is modified \citep{Jupp1975,CoffeyDepritDeprit1994}, or when other effects as lunisolar perturbations are included in the model \citep{Hough1981}. Also, in some cases the effects of higher order harmonics may be comparable to those of $J_2$, as happens for the selenopotential or the gravitational field of Venus. In these cases the disturbing effects of $J_2$ are not dominant, and critical inclinations, as well as the general orbit behavior, must be obviously discussed in the presence of all the first order effects \citep{LaraFerrerSaedeleer2009,Lara2011,LiuBaoyinMa2011}. In the present research the main problem of AST is the unique model under consideration and, consequently, the conclusions are limited to those cases in which the model applies. The subject of critical inclinations in different astrodynamical or astronomical contexts is not discussed here, and interested readers may consult the updated collection of references provided in the Introduction of \citep{RahomaKhattabAbdElSalam2014}.

\par

Aerospace engineers very soon found practical applications for satellite missions at the critical inclination, like the well-known Molniya or Tundra orbits \citep[see][for instance]{StoneBrodsky1988,Ulybyshev2009,BarkerStoen2001,ZhuZhaoZhang2014}. But useful orbits at the critical inclination are not limited to the case of highly eccentric orbits, and the critical inclination has been recently identified as a suitable choice for deploying cluster missions requiring bounded satellite motion \citep{LaraGurfil2012}.
\par

From a mathematical perspective, the critical inclination problem is commonly presented as a singularity in the solution of the motion of a massless particle in the gravitational field of an oblate body. This singularity is caused by the appearance of small divisors in the analytic integration of the secular terms of a perturbation theory developed by averaging \citep[see Sec.~12 of Chap.~17 of][for instance]{BrouwerClemence1961}. Eventually, some controversy arose on wether the singularity found at the critical inclination was just virtual or not, which was prolonged for some time \citep{Brouwer1963,Jupp1988}. Thus, on the one hand, the problem of small divisors is related with the occurrence of mean motion resonances between different frequencies of the motion; this connection is fairly clear in the case of orbital or tesseral resonances \citep{Garfinkel1982,FerrazMello1988}, but the frequencies involved in the resonance event at the critical inclination were not so apparent from Brouwer's action-angles approach to the analytic integration of the main problem. On the other hand, practitioners had not found the expected increase in the coordinates perturbations of orbits close to the critical inclination at those times \citep{Lubowe1969}. But the critical inclination singularity is, undoubtedly, essential. Indeed, following the lines of global geometric solutions proposed by \cite{Deprit1983,Deprit1984}, it was demonstrated by \cite{Cushman1983,Cushman1984} and \cite{CoffeyDepritMiller1986}, with later amendments by \cite{FerrerSanJuanAbad2007}, that the critical inclination phenomenon is produced by a change in the stability of circular orbits in a bifurcation event, a result which is in agreement with the behavior that had been anticipated by \cite{Izsak1962}.
\par

Remark that the troubles in integrating orbits at, or close to, the critical inclination happen only in the analytical approach, whereas there are no difficulties in the numerical or semi-analytical integration of orbits at the critical inclination. Furthermore, the small denominators problem is avoided in practical implementations of Brouwer's-type analytical integration by using special ways of handling the critical terms, cf.~Section 7 of \citep{CoffeyNealSegermanTravisano1995} or appendix A.F.~of \citep{HootsSchumacherGlover2004}. Even these days, analogous methods are proposed to cope with virtual singularities, like the cases of small eccentricities or inclinations, and the critical inclination \citep{XuXu2013}.
\par

Most mentioned research efforts base on averaging procedures, a fact that may prompt the belief that the critical inclination singularity only discloses in the treatment of the secular terms of the gravitational potential. However, the bifurcation event is found to become apparent too in the direct numerical integration of the main problem including both short- and long-period effects \citep{Broucke1994}.
Furthermore, the occurrence of critical inclinations in the analytical integration of spherical-variables intermediaries of the main problem is also well-known \citep{Sterne1958,Garfinkel1958,Aksnes1965}. For the latter, the critical inclination was clearly identified with a 1 to 1 commensurability between the satellite's draconitic (from ascending node to ascending node) and anomalistic (from perigee to perigee) periods, in this way providing a clear physical explanation of the resonance phenomenon. Namely, at resonance the sub-satellite point always reaches a given latitude in the same time interval (some multiple of the nodal period), and each time this happens the satellite's altitude over the surface of the earth is exactly the same. Therefore, because of the axial symmetry of the main problem potential, the satellite undergoes an identical gravitational pull after a constant period.
\par

Typpical intermediary orbits of the main problem can be integrated in closed form, but at the expense of using elliptic integrals \citep{Sterne1958,Garfinkel1958,Aksnes1965,CidLahulla1969}. Therefore, critical inclinations are only uncovered after a series expansion of corresponding closed form solutions. On the contrary, Deprit's radial intermediary \citep{Deprit1981} accepts a closed form solution which is free from elliptic integrals. This solution is used here to show that orbit inclination can be written explicitly as a function of the frequencies of the motion without need of resorting to averaging or series expansions. In this way, it is possible to demonstrate that resonances between the periods with which the polar and the radial variables evolve in the instantaneous orbital plane are parametrized by inclination. Besides, for the small values of $J_2$ which are characteristic of AST, it is shown that the 1 to 1 resonance is the unique deep resonance that can be guaranteed to exist. Finally, the existence of other critical inclinations is illustrated, which may happen only if much higher values of $J_2$ are considered.
\par

The paper is organized as follows. First, basic facts of the mean elements approach to the critical inclination resonance, limited to the first order of $J_2$, are summarized. Based on Battin's exposition, it is shown that inclination resonances may exist depending on the mean rate of rotation of the line of apsides. Then, after disclosing the relevance of using polar-nodal variables in the description of this type of resonance, the solution of Deprit's radial intermediary is reworked, but only for the trajectory in the instantaneous plane of motion. This extremely simple solution is finally used to find the closed form relation between inclination and the frequencies of the orbital motion in non-averaged, polar variables.

\section{Mean Elements Approach to the Critical Inclination}

Artificial satellite theory studies the motion of a massless particle in the presence of the gravitational potential. For earth-like bodies, the second order zonal harmonic coefficient ($J_2$) dominates all other harmonic coefficients, and hence the truncation of the expansion of the gravitational potential by neglecting all harmonic coefficients except $J_2$ is customarily called the \emph{main problem} of AST \citep{BrouwerClemence1961}.
\par

In Hamiltonian form, the main problem Hamiltonian is written
\begin{equation} \label{HamCartesian}
\mathcal{H}=\frac{1}{2}(X^2+Y^2+Z^2)-\frac{\mu}{r}+J_2\,\frac{\mu}{r}\,\frac{\alpha^2}{r^2}\,P_2(\sin\varphi),
\end{equation}
where the gravitational parameter $\mu$, the scaling factor $\alpha$ (the equatorial radius of the earth), and the oblateness coefficient $J_2$ are physical parameters that define the gravity field, $(x,y,z)$ are Cartesian coordinates, $(X,Y,Z)$ their conjugate momenta, respectively, which coincide with the Cartesian coordinates of the velocity in an inertial frame, $r=\sqrt{x^2+y^2+z^2}$ is the distance from the origin of coordinates, $\sin\varphi=z/r$, and $P_2$ is the Legendre polynomial of degree 2. In AST the oblateness coefficient is a ``small parameter'' which for the earth is of the order of one thousandth.
\par

The main problem is known to be non-integrable \citep{Danby1968,IrigoyenSimo1993,CellettiNegrini1995}, but a lot of information on the dynamics can be obtained under simplifying assumptions. The classical approach is based on expanding the disturbing function as a Fourier series on the mean anomaly, and retaining only the constant term \citep[cf.~Sec.~10.6 of][]{Battin1999}. Then, the mean effects on the orbital elements are integrated with the method of variations of parameters. In this way it is easily seen that the Lagrange planetary equation for the semi-major axis $a$, eccentricity $e$, and inclination $i$ vanish, whereas the variations of the other (mean) orbital elements can be integrated by quadrature.
\par

In particular, from Eq.~(10.95) in Sec.~10.6 of Battin's book \citep{Battin1999} we find
\begin{equation} \label{nwBattin}
\frac{{n}_\omega}{n}=\frac{3}{4}\sigma\,(5\cos^2i-1),
\end{equation}
where ${n}_\omega$ is the mean rate of rotation of the line of apsides, $n$ is the mean motion, %mean "mean motion" $n=\sqrt{\mu/a}$
\begin{equation} \label{sigma}
\sigma=J_2\,\frac{\alpha^2}{p^2},
\end{equation}
and
\begin{equation}
p=a\,(1-e^2),
\end{equation}
is the parameter of the conic or \emph{semilatus rectum}. Equation (\ref{nwBattin}) shows the existence of a critical inclination $i_\mathrm{c}$ such that $\cos^2i_\mathrm{c}=1/5$, at which ${n}_\omega=0$ and hence the line of apsides remains ``frozen'' on average. Otherwise the perigee regresses or advances depending on the sign of the coefficient $(5\cos^2i-1)$ in Eq.~(\ref{nwBattin}).
\par

\begin{figure}[htbp] \centering
\includegraphics[scale=1.3]{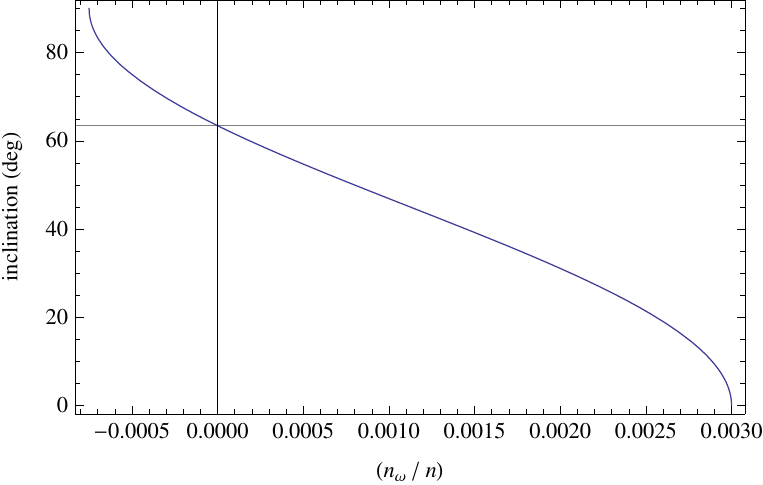}
\caption{Inclinations resonances according to Eq.~(\ref{inwn}) (low earth orbits).}
\label{f:icleo}
\end{figure}

Alternatively, Eq.~(\ref{nwBattin}) can be solved for the inclination, giving
\begin{equation} \label{inwn}
i=\arccos\sqrt{\frac{1}{5}\left(1+\frac{4}{3}\frac{1}{\sigma}\frac{{n}_\omega}{n}\right)},
\end{equation}
which shows that orbital inclination can be written explicitly as a function of the ratio of the mean rate of rotation of the line of apsides to the orbital mean motion. Specifically, Eq.~(\ref{inwn}) shows that resonances between the mean rate of rotation of the perigee and the mean motion are very shallow for the artificial satellite problem, and hence are not of concern in AST. Indeed, for low earth orbits $\sigma\sim{J}_2=\mathcal{O}(10^{-3})$ and the resonance condition will occur only after hundreds of satellite orbits in the more favorable cases, as illustrated in Fig.~\ref{f:icleo}.

\begin{figure}[htbp] \centering
\includegraphics[scale=1.3]{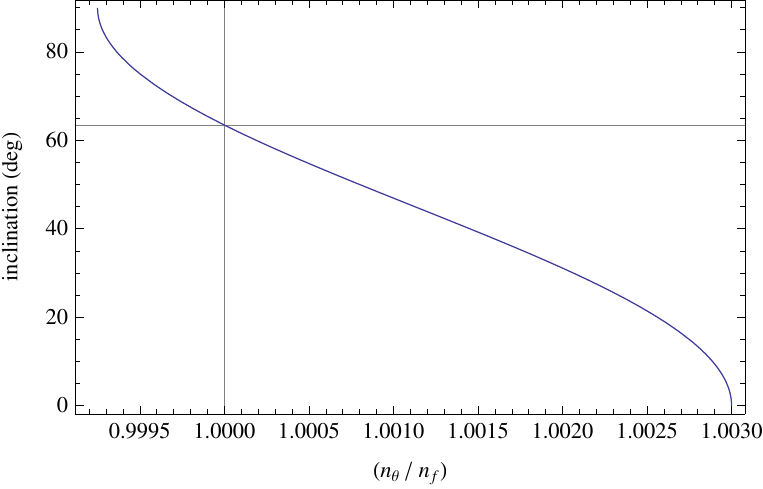}
\caption{The critical inclination as an inclinations resonance (low earth orbits).}
\label{f:inznf}
\end{figure}

Still, in spite of the clear physical interpretation of frozen-perigee orbits as resonant orbits, if the orbital element set chosen for representing the mean elements solution includes the argument of the perigee among its variables\footnote{A survey of the most common  sets of orbital elements can be found in \citep{Hintz2008}} then the critical inclination does not manifest as the usual resonance between different frequencies of the motion. However, it is readily seen that the orbital frequencies that are involved in inclination resonances are clearly apparent by simply recalling that
\begin{equation} \label{wzf}
\omega=\theta-f,
\end{equation}
where, $\theta$ is the argument of latitude and $f$ the true anomaly. Indeed, in view of the average rates with which the true and mean anomaly evolve are the same, Eq.~(\ref{inwn}) can be rewritten as
\begin{equation}\label{inzn}
i=\arccos\sqrt{\frac{1}{5}\left[1+\frac{4}{3}\frac{1}{\sigma}\left(\frac{{n}_\theta}{{n}_f}-1\right)\right]},
\end{equation}
where ${n}_\theta$ is the mean rate of rotation of the argument of the latitude and ${n}_f$ is the mean rate of rotation of the true anomaly. Then, the diagram of Fig.~\ref{f:icleo} is rendered again but now from Eq.~(\ref{inzn}), in this way displacing abscissas to the right by one unit to obtain an inclination resonances diagram which now does include the critical inclination of $63.4\deg$ as the 1 to 1 resonance, as shown in Fig.~\ref{f:inznf}.
\par

Hence, the simple replacement provided by Eq.~(\ref{wzf}) suggests that the set of polar-nodal variables is a suitable orbital elements representation for illustrating the inclination resonance phenomenon.
\par

\section{Inclination Resonances without Averaging}

Because of the axial symmetry of the disturbing function (the term factored by the small parameter $J_2$), the main problem Hamiltonian is only of two degrees of freedom, a fact that becomes evident when formulating Eq.~(\ref{HamCartesian}) in the set of polar-nodal variables $(r,\theta,\nu,R,\Theta,N)$, which stand for the radius from the earth's center of mass, the argument of latitude, the Right Ascension of the ascending node, the radial velocity, the modulus of the angular momentum vector, and the polar component of the angular momentum vector, respectively. Then, Eq.~(\ref{HamCartesian}) is written
\begin{equation} \label{Ham}
\mathcal{H}=\frac{1}{2}\left(R^2+\frac{\Theta^2}{r^2}\right)-\frac{\mu}{r}\left[1
+J_2\,\frac{\alpha^2}{r^2}\left(\frac{1}{2}-\frac{3}{4}\sin^2i+\frac{3}{4}\sin^2i\cos2\theta\right)\right],
\end{equation}
where, in the polar-nodal set of canonical variables, the orbit inclination is obtained from $\cos{i}=N/\Theta$.
\par

The fact that $\nu$ is cyclic in Eq.~(\ref{Ham}) reflects the axial symmetry accepted by the main problem. Indeed, its conjugate momentum $N$ is an integral, and, therefore, the reduced flow in the instantaneous plane of motion
\begin{equation} \label{flow2D}
\frac{\mathrm{d}(r,\theta)}{\mathrm{d}t}=\frac{\partial\mathcal{H}}{\partial(R,\Theta)},
\qquad
\frac{\mathrm{d}(R,\Theta)}{\mathrm{d}t}=-\frac{\partial\mathcal{H}}{\partial(r,\theta)},
\end{equation}
which is of two degrees of freedom, decouples from the rotation of the node
\begin{equation} \label{noderate}
\frac{\mathrm{d}\nu}{\mathrm{d}t}=\frac{\partial\mathcal{H}}{\partial{N}},
\end{equation}
which can be integrated by quadrature after solving the differential system in polar variables given in Eq.~(\ref{flow2D}).
\par

In spite of the reduced problem remains non-integrable, since it is only of two degrees of freedom the computation of Poincar\'e surfaces of section and periodic orbits may be used to get insight in the dynamics of Eq.~(\ref{flow2D}) \citep{Danby1968,Broucke1994}. Alternatively, the desired information on the dynamics may be obtained by computing approximated solutions to the flow.
\par

\subsection{The Radial Intermediary}

The use of a variety of integrable approximations of the main problem, the so-called \emph{intermediaries}, have been proposed in the literature for approximating the main problem dynamics. In particular, this work deals exclusively with Deprit's natural, radial intermediary
\begin{equation} \label{HamP}
\mathcal{H}=\frac{1}{2}\left(R^2+\frac{Q^2}{r^2}\right)-\frac{\mu}{r},
\end{equation}
where $Q\equiv{Q}(\Theta,N)$ is a constant function given by
\begin{equation} \label{Q}
Q=\Theta\,\sqrt{1+\sigma\left(\frac{1}{2}-\frac{3}{2}\cos^2i\right)},
\end{equation}
and $\sigma$ was defined in Eq.~(\ref{sigma}) where, now, the \textit{semilatus rectum}
\begin{equation} \label{ppn}
p=\frac{\Theta^2}{\mu},
\end{equation}
is expressed in polar-nodal variables.
\par
 
Equation (\ref{HamP}) is obtained after applying the elimination of the parallax transformation to Eq.~(\ref{Ham}). Up to the first order of $J_2$, it converts the main problem into a quasi-Keplerian, integrable problem \citep{Deprit1981}. For the reader's convenience, the contact transformation that converts the main problem Hamiltonian into Deprit's radial intermediary in Eq.~(\ref{HamP}) is provided in the appendix. Note that the expressions provided in the appendix are simpler than those given in Deprit's original paper. Indeed, in view of the recent claims that the elimination of the parallax transformation is naturally conceived, as well as more easily accomplished, when working in the setting provided by Delaunay variables \citep{LaraSanJuanLopezOchoa2013b,LaraSanJuanLopezOchoa2013c}, the use of the $C$ and $S$ parallactic functions, which are essential in the Deprit's construction of the algorithm in polar-nodal varaibles, does not add any value to this canonical transformation. Besides, the expressions provided in the appendix are better tuned for fast evaluation which may be an essential prerequisite in practical applications \citep{GurfilLara2014}.
\par

A brief outline of the conventional integration of the Hamiltonian flow of Eq.~(\ref{HamP}) is given in \citep{Deprit1981}, where the alternative integration by means of a \emph{torsion} transformation is also provided with much more detail. The latter being a transformation in implicit variables, Deprit reconstructed the torsion as an explicit transformation by means of a Lie transform which, up to the first order of $J_2$, shows that the torsion in the orbital plane becomes the identity mapping at the critical inclination. Here we take a different approach in which, without need of resorting to series expansions, the trajectory solution in polar variables is used to demonstrate that orbit inclination can be expressed as a function of the ratio between the draconitic and anomalistic periods or frequencies. This closed form relation is then used to compute the (critical) inclination at which both periods get the same value, therefore meeting the frozen orbit condition.
\par

\subsection{Trajectory in the Instantaneous Plane of Motion}

Because the argument of the latitude is cyclic in Eq.~(\ref{HamP}) the modulus of the angular momentum is constant, and, consequently, the orbital plane evolves with constant inclination in the transformed phase space. Then, the one degree of freedom problem in $(r,R)$ decouples from the motion of $\theta$ and $\nu$.
The reduced problem is easily integrated by noting that, from Hamilton equations,
\begin{equation} \label{rR}
\frac{\mathrm{d}r}{\mathrm{d}t}=\frac{\partial\mathcal{H}}{\partial{R}}=R.
\end{equation}
Then, for any manifold $\mathcal{H}=h$, Eq.~(\ref{rR}) is replaced into Eq.~(\ref{HamP}) leading to
\begin{equation} \label{rp}
\frac{\mathrm{d}r}{\mathrm{d}t}=\sqrt{2h+2\frac{\mu}{r}-\frac{Q^2}{r^2}}.
\end{equation}
\par

Equation (\ref{rp}) is in separate variables and the independent variable is easily integrated by quadrature using the standard change of variable from the radius to the eccentric anomaly.
Besides,
\begin{equation} \label{zp}
\frac{\mathrm{d}\theta}{\mathrm{d}t}=\frac{\partial\mathcal{H}}{\partial\Theta}=\frac{P}{r^2},
\end{equation}
where $P\equiv{P}(\Theta,N)$ is a constant function given by
\begin{equation} \label{P}
P=\Theta\left[1-\sigma\left(\frac{1}{2}-3\cos^2i\right)\right].
\end{equation}
\par

Therefore, the equation of the trajectory in the instantaneous plane of motion is obtained combining Eqs.~(\ref{zp}) and (\ref{rp}) into
\begin{equation} \label{zr}
\frac{\mathrm{d}\theta}{\mathrm{d}r}=
\frac{P}{r^2\,\sqrt{2h+2(\mu/r)-(Q^2/r^2)}},
\end{equation}
which is also in separate variables and hence can be integrated by quadrature to give
\begin{equation} \label{z}
\theta=P\,\int_{r_{\mathrm{m}}}^r\frac{\mathrm{d}s}{ s^2\sqrt{ 2h+2(\mu/s)-(Q^2/s^2) } },
\end{equation}
where $r_{\mathrm{m}}$ is the minimum value of the radial distance.
\par

The trajectory in the orbital plane given in Eq.~(\ref{z}) is solved by the standard change of variable from the radius to the true anomaly $f$. Indeed, making
\begin{equation} \label{ellipse}
r=\frac{Q^2/\mu}{1+e\cos{f}},
\end{equation}
where
\begin{equation} \label{eccentricity}
e^2=1+\frac{2h}{\mu}\,\frac{Q^2}{\mu},
\qquad
0\le{e}\le1,
\end{equation}
and taking into account that $r=r_\mathrm{m}\Rightarrow{f}=0$
, one trivially gets
\begin{equation} \label{zeta}
\theta=\theta_0+\frac{P}{Q}\,f.
\end{equation}
Hence, using Eqs.~(\ref{ellipse}), (\ref{zeta}), and (\ref{ppn}), the trajectory in the orbital plane is expressed as
\begin{equation} \label{trajectory}
r=\frac{p\,(Q/\Theta)^2}{1+e\cos[(Q/P)(\theta-\theta_0)]}.
\end{equation}
\par

\subsection{Critical inclinations}

By differentiation of Eq.~(\ref{zeta})
\begin{equation} \label{res}
\frac{n_r}{n_\theta}=\frac{Q}{P}.
\end{equation}
where $n_\theta=\mathrm{d}\theta/\mathrm{d}t$ is the draconitic frequency, that is $n_\theta=2\pi/T_\theta$ where $T_\theta$ is the time elapsed between two consecutive passages of the satellite through the ascending node, and $n_r=\mathrm{d}f/\mathrm{d}t$ is the anomalistic frequency, that is $n_r=2\pi/T_r$ where $T_r$ is the time elapsed between two consecutive satellite's transits through the perigee.
\par

\begin{figure}[htbp]\centering %\centerline{
\includegraphics[scale=1.3]{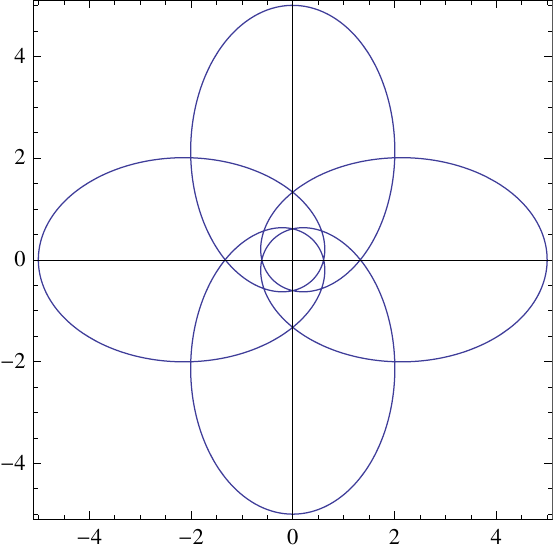} \\[1ex]
\includegraphics[scale=1.3]{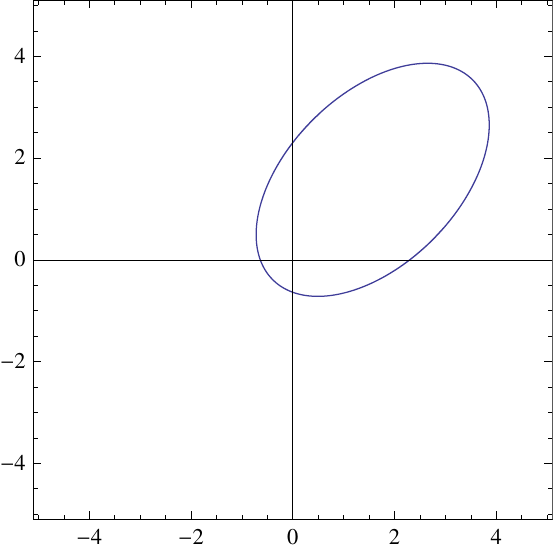} \\[1ex]
\includegraphics[scale=1.3]{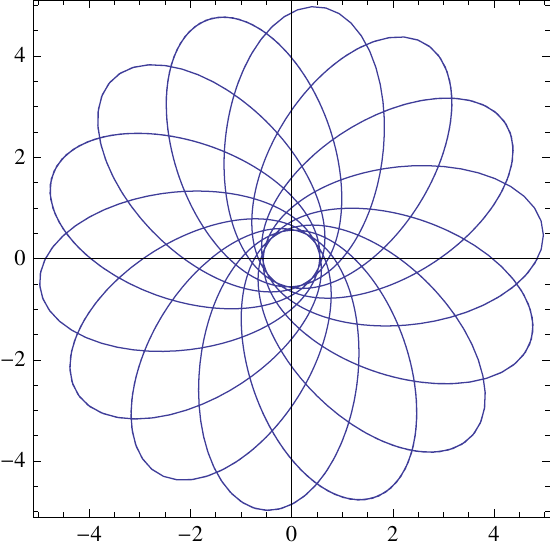}
%}
\caption{Sample solutions of Eq.~(\protect\ref{trajectory}) for $Q/P=4/5$ (top), $1$ (center), and $14/13$ (bottom).}
\label{f:rosettes}
\end{figure}

In general, $P$ and $Q$ will be incommensurable numbers and the trajectory in the orbital plane given by Eq.~(\ref{trajectory}) will be a rosette with an infinite number of perigees. However, the frequencies $n_\theta$ and $n_r$ will become commensurable for those values of $\Theta$ and $N$ that convert Eq.~(\ref{res}) into a rational number $k\equiv{k}(\Theta,N)$. When this happens the rosette becomes a periodic orbit which closes after $\theta$ has walked a number $P$ of cycles ---as illustrated in Fig.~\ref{f:rosettes} for, from top to bottom, $Q/P=4/5$, $1$, and $14/13$ ($e=0.8$, $Q^2=\mu$, $\theta_0=3\pi/4$).
\par

If $Q/P$ in the right hand of Eq.~(\ref{res}) is now replaced by the division of Eq.~(\ref{Q}) by Eq.~(\ref{P}), it results
\begin{equation} \label{ratio}
\frac{n_r}{n_\theta}=\frac{\sqrt{1+\sigma\left(\frac{1}{2}-\frac{3}{2}\cos^2i\right)}}{1-\sigma\left(\frac{1}{2}-3\cos^2i\right)},
\end{equation}
from which orbit inclination can be written explicitly as a function of the frequencies of the reduced problem. Namely,
\begin{equation} \label{cos2i}
\cos^2i=\frac{\sqrt{1+4(6+\sigma)\,(n_r/n_\theta)^2}-1-2(2-\sigma)\,(n_r/n_\theta)^2}{12\sigma\,(n_r/n_\theta)^2}.
\end{equation}
\par

Leaving aside the case of almost rectilinear orbits, following results limit to non-impact orbits, from which $p\ge\alpha$.
In these cases $\sigma$ is small, and arbitrary rational values $(n_r/n_\theta)=k$ will not result, in general, in real inclinations producing periodicity in the instantaneous plane of motion. Therefore, one must carefully explore the range of rational values of $k$ which are allowed in Eq.~(\ref{cos2i}). 
\par

Expansion of Eq.~(\ref{cos2i}) in power series of $\sigma$ gives
\begin{equation} \label{cos2is}
\cos^2i=\frac{\sqrt{1+24k^2}-1-4 k^2}{12k^2}\,\sigma^{-1}
+\frac{1}{6}\left(1+\frac{1}{\sqrt{1+24k^2}}\right)
-\frac{k^2}{6(1+24k^2)^{3/2}}\,\sigma+\mathcal{O}(\sigma^2)
%+\frac{\sigma ^2 k^4}{3 \left(24k^2+1\right)^{5/2}}+O\left(\sigma ^3\right)
\end{equation}
which clearly shows that $k^2=1$, corresponding to the 1 to 1 resonance between the draconitic and anomalistic periods, is the only value of $k$ that cancels the coefficient of $\sigma^{-1}$ ---the first term in the right side of Eq.~(\ref{cos2is}). 
\par

Therefore, the closed form expression in Eq.~(\ref{cos2i}) is simplified for the 1 to 1 resonance to
\begin{equation} \label{cos2i1to1}
\cos^2i_\mathrm{c} = \frac{1}{6}-\frac{5}{12\sigma}\left(1-\sqrt{1+\frac{4}{25}\sigma}\right),
\end{equation}
which shows the dependence of the critical inclination on $\sigma$ and, through it, on both the $J_2$ value and the modulus of the angular momentum, cf.~Eqs.~(\ref{sigma}) and (\ref{ppn}). The limit $\sigma\rightarrow0$ results in the well-known value of the critical inclination $\cos^2i_\mathrm{c}=1/5$. For the earth, in which case $J_2=\mathcal{O}(10^{-3})$, it is found that this critical inclination value is accurate to the order of $J_2^2$. 
Indeed,
\begin{equation} \label{cos2i1to1s}
\cos^2i_\mathrm{c}=\frac{1}{5}-\frac{1}{750}\sigma+\frac{1}{9375}\sigma^2+\mathcal{O}(\sigma^3),
\end{equation}
where $(1/750)\sim{J}_2$, and $J_2\sim\sigma$ for low earth orbits.
\par

Besides, because $\theta=\omega+f$, Eq.~(\ref{zeta}) trivially shows that in the 1 to 1 resonance the perigee becomes fixed (or ``frozen'') at $\omega=\theta_0$. The degeneracy of the solution (all orbits become frozen at the critical inclination) is a result of the integrability of the Hamiltonian truncation to the first order of $J_2$ in Eq.~(\ref{HamP}). However, as it is well-known, second order effects of $J_2$ break this degeneracy to leave only a discrete number of frozen orbits \citep[see][and references therein]{CoffeyDepritMiller1986}.
\par

Other possible resonances would also require $k^2\approx1$ to compensate the smallness of $\sigma$ in the first summand of Eq.~(\ref{cos2is}) so that the condition $0\le\cos^2i\le1$ can be fullfiled. Because this only might happen for higher order resonances $k\sim1+\mathcal{O}(\sigma)$, corresponding inclinations, if they exist, will not cause any small-divisor type complication in a perturbation theory, and hence are not of major concern in AST.
\par

\begin{figure}[htbp] \centering
\includegraphics[scale=1.5]{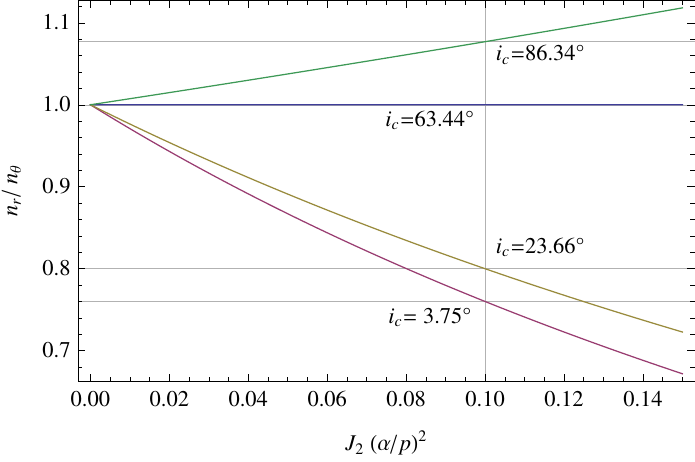}
\caption{Evolution of $k=(n_r/n_\theta)$ with $\sigma=J_2\,(\alpha/p)^2$ for different inclinations, from Eq.~(\protect\ref{ratio}).}
\label{f:ksigma}
\end{figure}

Conversely, assumed that higher values of $J_2$ may exist, other critical inclinations can be found. Thus, for instance, for $\sigma=0.1$ ($J_2\sim\mathcal{O}(0.1)$) the resonances $n_r/n_\theta=19/25$, $4/5$, $1/1$, and $14/13$, will result in critical inclinations at $3.75$, $23.66$, $63.43$, and $86.34$ deg, respectively, as obtained from Eq.~(\ref{cos2i}). This is illustrated in Fig.~\ref{f:ksigma}, where the ratio $n_r/n_\theta$, as given in Eq.~(\ref{ratio}), is displayed as function of $\sigma$ for different inclinations. Ordinates $n_r/n_\theta=0.76$, $0.8$, $1$, and $1.077$ correspond to the previously mentioned resonances for the abscissa $\sigma=0.1$. Remark that none of the inclination lines in Fig.~\ref{f:ksigma} is a straight line, not even the line of $63.44\deg$.
\par

\section{Conclusions}

The mapping between orbital inclination and the ratio draconitic anomalistic frequencies of the main problem of AST can be made explicit in both directions in closed form without need of relying on averaging procedures. This fact lets look at the small divisors problem arising in the vicinity of the critical inclination as the familiar problem of resonances between the different frequencies of the motion, which in this particular case happen between the rates of variation of both polar motion variables.
For those inclinations leading to resonance the orbit in the instantaneous plane of motion turns into a closed rosette. In the particular case of the 1 to 1 resonance the trajectory in the orbital plane becomes a mere ellipse which, therefore, has the perigee frozen in spite of its non-Keplerian character.
\par

Analogously to the case of tesseral resonances, for which resonant motion is parametrized by semi-major axis, resonant satellite motion in axial-symmetric force fields is parametrized by inclination. For the latter, the only relevant resonance in artificial satellite theory is the 1 to 1 resonance. In the case of the earth, this resonance happens at the well-known value of the critical inclination, which is accurate up to second-order effects of $J_2$. However, there are no mathematical objections for other critical inclinations of the main problem to occur assumed that the $J_2$ coefficient may take much higher values than those of earth-like bodies.

\section*{Acknowledgemnts}

Part of this research has been supported by the Government of Spain (Project AYA 2010-18796).

\appendix

\section{Transformation to the Radial Intermediary}

Equation (\ref{HamP}) is obtained after the elimination of the parallax ---a canonical transformation from original polar variables $(r,\theta,\nu,R,\Theta,N)$ to new variables $(r',\theta',\nu',R',\Theta',N')$--- is carried out up to the first order of $J_2$. Therefore,  Eq.~(\ref{HamP}) should be written in prime variables.
\par

Corresponding transformation equations are
\begin{eqnarray} \label{Dr}
\frac{r-r'}{p'} &=& \kappa\left(1-\frac{3}{2}s'^2-\frac{1}{2}s'^2\cos2\theta'\right) \\[1ex] %\label{Dz}
\theta-\theta' &=& \kappa\left\{
\left[\frac{3}{4}-\frac{5}{4}c'^2-(1-3c'^2)\frac{p'}{r'}\right]\!\sin2\theta' 
+\frac{p'R'}{\Theta'}\left[1-6c'^2+(1-2c'^2)\cos2\theta'\right]\!\right\} \qquad \\[1ex] %\label{Dv}
\nu-\nu' &=& \kappa\,c'\left[
 \left(\frac{1}{2}-2\,\frac{p'}{r'}\right)\!\sin2\theta'+\frac{p'R'}{\Theta}\,(3+\cos2\theta') \right] \\[1ex] %\label{DR}
\frac{R-R'}{\Theta'/p'} &=& \kappa\,\frac{p'^2}{r'^2}\,s'^2\sin2\theta' \\[1ex] %\label{DZ}
\frac{\Theta-\Theta'}{\Theta'} &=& \kappa\,s'^2\!\left[ \left(\frac{1}{2}-2\,\frac{p'}{r'}\right)\!\cos2\theta'-\frac{p'R'}{\Theta'}\sin2\theta' \right] \\[1ex]  \label{DN}
N &=& N'
\end{eqnarray}
where
\begin{equation}
\kappa%=-\frac{1}{2}\sigma
=-\frac{1}{2}J_2\,\frac{\alpha^2}{p'^2},
\qquad
p'=\frac{\Theta'^2}{\mu},
\end{equation}
and $s'\equiv\sin{i}'$, $c'\equiv\cos{i}'=N/\Theta'$.
\par

Note that, instead of taking the original transformation equations from \citep{Deprit1981} (see p.~133), we borrowed the transformation equations from \citep{GurfilLara2014} because of their simplicity: they only require the evaluation of sine and cosine functions of argument $2\theta$, whereas trigonometric functions of $\theta$, $2\theta$ and $3\theta$ need to be evaluated when using the transformation equations in the original form given by Deprit. These simpler expressions were rendered after realizing that the so-called parallactic functions $C$ and $S$, two invariants of the Keplerian motion pertaining to the kernel of the Lie derivative in the algebra of functions in which the original elimination of the parallax is based, are not necessary for carrying out the elimination of the parallax. Indeed, straightforward derivations demonstrate the higher efficiency of the elimination of the parallax when approached in Delaunay variables \citep{LaraSanJuanLopezOchoa2013b,LaraSanJuanLopezOchoa2013c}. Then, after computing the elimination of the parallax in the Delaunay chart, the generating function of this Hamiltonian \emph{simplification} \citep{DepritFerrer1989} is easily reformulated in the more convenient set of polar-nodal variables by means of standard relations, from which the transformation equations (\ref{Dr})--(\ref{DN}) are readily obtained.

%\bibliographystyle{/Users/martinlara/Research/bibrefs/elsarticle-harv}   % Number the references.
%%\bibliographystyle{/Users/martinlara/Research/bibrefs/elsarticle-num}   % Number the references.
%%\footnotesize
%\bibliography{/Users/martinlara/Research/bibrefs/references}   % Use references.bib to resolve the labels.

\section*{References}
\begin{thebibliography}{45}
\expandafter\ifx\csname natexlab\endcsname\relax\def\natexlab#1{#1}\fi
\expandafter\ifx\csname url\endcsname\relax
  \def\url#1{\texttt{#1}}\fi
\expandafter\ifx\csname urlprefix\endcsname\relax\def\urlprefix{URL }\fi

\bibitem[{{Aksnes}(1965)}]{Aksnes1965}
{Aksnes}, K., Aug. 1965. {On the dynamical theory of a near-earth satellite,
  I.} Astrophysica Norvegica 10, 69--77.

\bibitem[{Barker and Stoen(2001)}]{BarkerStoen2001}
Barker, L., Stoen, J., 2001. {Sirius satellite design: the challenges of the
  Tundra orbit in commercial spacecraft design}. In: Culp, R.~D., Schira, C.~N.
  (Eds.), Guidance and Control 2001. Vol. 107 of Advances in the Astronautical
  Sciences. American Astronautical Society, Univelt, Inc., USA, pp. 575--596.

\bibitem[{Battin(1999)}]{Battin1999}
Battin, R.~H., 1999. An Introduction to the Mathematics and Methods of
  Astrodynamics. American Institute of Aeronautics and Astronautics, Reston,
  VA.

\bibitem[{{Broucke}(1994)}]{Broucke1994}
{Broucke}, R.~A., Feb. 1994. {Numerical integration of periodic orbits in the
  main problem of artificial satellite theory}. Celestial Mechanics and
  Dynamical Astronomy 58~(2), 99--123.

\bibitem[{{Brouwer}(1963)}]{Brouwer1963}
{Brouwer}, D., 1963. {Review of Celestial Mechanics}. Annual Review of
  Astronomy and Astrophysics 1, 219--234.

\bibitem[{Brouwer and Clemence(1961)}]{BrouwerClemence1961}
Brouwer, D., Clemence, G.~M., 1961. Methods of Celestial Mechanics. Academic
  Press, New York and London.

\bibitem[{{Celletti} and Negrini(1995)}]{CellettiNegrini1995}
{Celletti}, A., Negrini, P., March 1995. Non-integrability of the problem of
  motion around an oblate planet. Celestial Mechanics and Dynamical Astronomy
  61, 253--260.

\bibitem[{Cid and Lahulla(1969)}]{CidLahulla1969}
Cid, R., Lahulla, J.~F., 1969. {Perturbaciones de corto periodo en el
  movimiento de un sat\'{e}lite artificial, en funci\'on de las variables de
  Hill}. Publicaciones de la Revista de la Academia de Ciencias de Zaragoza 24,
  159--165.

\bibitem[{{Coffey} et~al.(1994){Coffey}, {Deprit}, and
  {Deprit}}]{CoffeyDepritDeprit1994}
{Coffey}, S.~L., {Deprit}, A., {Deprit}, E., May 1994. {Frozen orbits for
  satellites close to an earth-like planet}. Celestial Mechanics and Dynamical
  Astronomy 59~(1), 37--72.

\bibitem[{{Coffey} et~al.(1986){Coffey}, {Deprit}, and
  {Miller}}]{CoffeyDepritMiller1986}
{Coffey}, S.~L., {Deprit}, A., {Miller}, B.~R., Dec. 1986. {The critical
  inclination in artificial satellite theory}. Celestial Mechanics 39~(4),
  365--406.

\bibitem[{Coffey et~al.(1996)Coffey, Neal, Segerman, and
  Travisano}]{CoffeyNealSegermanTravisano1995}
Coffey, S.~L., Neal, H.~L., Segerman, A.~M., Travisano, J.~J., 1996. An
  analytic orbit propagation program for satellite catalog maintenance. In:
  Alfriend, K.~T., Ross, I.~M., Misra, A.~K., Peters, C.~F. (Eds.), AAS/AIAA
  Astrodynamics Conference 1995. Vol.~90 of Advances in the Astronautical
  Sciences. American Astronautical Society, Univelt, Inc., USA, pp. 1869--1892.

\bibitem[{{Cushman}(1983)}]{Cushman1983}
{Cushman}, R., Dec. 1983. {Reduction, Brouwer's Hamiltonian, and the critical
  inclination}. Celestial Mechanics 31~(4), 401--429.

\bibitem[{{Cushman}(1984)}]{Cushman1984}
{Cushman}, R., Jul. 1984. {Erratum: ''Reduction, Brouwer's Hamiltonian, and the
  critical inclination'' [Celest. Mech., Vol. 31, No. 4, p. 401 - 429 (1983)].}
  Celestial Mechanics 33~(3), 297--297.

\bibitem[{{Danby}(1968)}]{Danby1968}
{Danby}, J.~M.~A., Dec. 1968. {Motion of a Satellite of a Very Oblate Planet}.
  The Astronomical Journal 73~(10), 1031--1038.

\bibitem[{Deprit(1981)}]{Deprit1981}
Deprit, A., 1981. The elimination of the parallax in satellite theory.
  Celestial Mechanics 24~(2), 111--153.

\bibitem[{{Deprit}(1983)}]{Deprit1983}
{Deprit}, A., Mar. 1983. {The reduction to the rotation for planar perturbed
  Keplerian systems}. Celestial Mechanics 29, 229--247.

\bibitem[{{Deprit}(1984)}]{Deprit1984}
{Deprit}, A., 1984. {Dynamics of orbiting dust under radiation pressure}. In:
  {Berger}, A. (Ed.), Big-Bang Cosmology Symposium in honour of G. Lemaitre.
  pp. 151--180.

\bibitem[{{Deprit} and {Ferrer}(1989)}]{DepritFerrer1989}
{Deprit}, A., {Ferrer}, S., Dec. 1989. {Simplifications in the theory of
  artificial satellites}. Journal of the Astronautical Sciences 37, 451--463.

\bibitem[{{Ferraz-Mello}(1988)}]{FerrazMello1988}
{Ferraz-Mello}, S., 1988. {On resonance}. Celestial Mechanics 43~(1-4), 69--89.

\bibitem[{{Ferrer} et~al.(2007){Ferrer}, {San-Juan}, and
  {Abad}}]{FerrerSanJuanAbad2007}
{Ferrer}, S., {San-Juan}, J.~F., {Abad}, A., Sep. 2007. {A note on lower bounds
  for relative equilibria in the main problem of artificial satellite theory}.
  Celestial Mechanics and Dynamical Astronomy 99~(1), 69--83.

\bibitem[{{Garfinkel}(1958)}]{Garfinkel1958}
{Garfinkel}, B., Mar. 1958. {On the motion of a satellite of an oblate planet}.
  The Astronomical Journal 63~(1257), 88--96.

\bibitem[{{Garfinkel}(1982)}]{Garfinkel1982}
{Garfinkel}, B., Nov. 1982. {On resonance in celestial mechanics (A survey)}.
  Celestial Mechanics 28~(3), 275--290.

\bibitem[{Gurfil and Lara(2014)}]{GurfilLara2014}
Gurfil, P., Lara, M., 2014. {Natural Intermediaries as Onboard Orbit
  Propagators (IAA-AAS-DyCoSS2-05-02)}. In: Second IAA Conference on Dynamics
  and Control of Space Systems, Rome, Italy, March 24-26, 2014.

\bibitem[{{Hagihara}(1961)}]{Hagihara1961}
{Hagihara}, Y., 1961. {Libration of an Earth Satellite with Critical
  Inclination}. Smithsonian Contributions to Astrophysics 5~(5), 39--51.

\bibitem[{Hintz(2008)}]{Hintz2008}
Hintz, G., May-June 2008. Survey of orbit element sets. Journal of Guidance,
  Control, and Dynamics 31~(3), 785--790.

\bibitem[{{Hoots} et~al.(2004){Hoots}, {Schumacher}, and
  {Glover}}]{HootsSchumacherGlover2004}
{Hoots}, F.~R., {Schumacher}, Jr., P.~W., {Glover}, R.~A., Mar. 2004. {History
  of Analytical Orbit Modeling in the U. S. Space Surveillance System}. Journal
  of Guidance, Control, and Dynamics 27~(5), 174--185.

\bibitem[{{Hough}(1981)}]{Hough1981}
{Hough}, M.~E., Oct. 1981. {Orbits near critical inclination, including
  lunisolar perturbations}. Celestial Mechanics 25~(2), 111--136.

\bibitem[{{Irigoyen} and {Simo}(1993)}]{IrigoyenSimo1993}
{Irigoyen}, M., {Simo}, C., March 1993. Nonintegrability of the {$J_2$}
  problem. Celestial Mechanics and Dynamical Astronomy 55~(3), 281--287.

\bibitem[{{Izsak}(1962)}]{Izsak1962}
{Izsak}, I.~G., Mar. 1962. {On the Critical Inclination in Satellite Theory}.
  SAO Special Report 90.

\bibitem[{{Jupp}(1975)}]{Jupp1975}
{Jupp}, A.~H., May 1975. {The problem of the critical inclination revisited}.
  Celestial Mechanics 11~(3), 361--378.

\bibitem[{{Jupp}(1988)}]{Jupp1988}
{Jupp}, A.~H., 1988. {The critical inclination problem - 30 years of progress}.
  Celestial Mechanics 43~(1-4), 127--138.

\bibitem[{{Kozai}(1961)}]{Kozai1961}
{Kozai}, Y., 1961. {Motion of a Particle with Critical Inclination in the
  Gravitational Field of a Spheroid}. Smithsonian Contributions to Astrophysics
  5~(5), 53--58.

\bibitem[{{Lara}(2011)}]{Lara2011}
{Lara}, M., Aug. 2011. {Design of long-lifetime lunar orbits: A hybrid
  approach}. Acta Astronautica 69~(3--4), 186--199.
%\newline\urlprefix\url{http://www.sciencedirect.com/science/article/pii/S00945%
%7651100066X}

\bibitem[{Lara et~al.(2009)Lara, Ferrer, and
  Saedeleer}]{LaraFerrerSaedeleer2009}
Lara, M., Ferrer, S., Saedeleer, B., 2009. Lunar analytical theory for polar
  orbits in a 50-degree zonal model plus third-body effect. The Journal of the
  Astronautical Sciences 57~(3), 561--577.
%\newline\urlprefix\url{http://dx.doi.org/10.1007/BF03321517}

\bibitem[{{Lara} and {Gurfil}(2012)}]{LaraGurfil2012}
{Lara}, M., {Gurfil}, P., Nov. 2012. {Integrable approximation of
  $J_{2}$-perturbed relative orbits}. Celestial Mechanics and Dynamical
  Astronomy 114~(3), 229--254.

\bibitem[{Lara et~al.(2014{\natexlab{a}})Lara, San-Juan, and
  L\'opez-Ochoa}]{LaraSanJuanLopezOchoa2013c}
Lara, M., San-Juan, J.~F., L\'opez-Ochoa, L.~M., 2014{\natexlab{a}}. {Delaunay
  variables approach to the elimination of the perigee in Artificial Satellite
  Theory}. Celestial Mechanics and Dynamical Astronomy in press~({also,
  arXiv:1312.7577v1 [nlin.CD]}).

\bibitem[{Lara et~al.(2014{\natexlab{b}})Lara, San-Juan, and
  L\'opez-Ochoa}]{LaraSanJuanLopezOchoa2013b}
Lara, M., San-Juan, J.~F., L\'opez-Ochoa, L.~M., 2014{\natexlab{b}}. {Proper
  Averaging Via Parallax Elimination (AAS 13-722)}. In: Astrodynamics 2013.
  Vol. 150 of {Advances in the Astronautical Sciences}. American Astronautical
  Society, Univelt, Inc., USA, pp. 315--331.

\bibitem[{{Liu} et~al.(2011){Liu}, {Baoyin}, and {Ma}}]{LiuBaoyinMa2011}
{Liu}, X., {Baoyin}, H., {Ma}, X., Jul. 2011. {Extension of the critical
  inclination}. Astrophysics and Space Science 334~(1), 115--124.
%\newline\urlprefix\url{http://dx.doi.org/10.1007/s10509-011-0685-y}

\bibitem[{{Lubowe}(1969)}]{Lubowe1969}
{Lubowe}, A.~G., Mar. 1969. {How Critical is the Critical Inclination?}
  Celestial Mechanics 1~(1), 6--10.

\bibitem[{{Rahoma} et~al.(2014){Rahoma}, {Khattab}, and {Abd
  El-Salam}}]{RahomaKhattabAbdElSalam2014}
{Rahoma}, W.~A., {Khattab}, E.~H., {Abd El-Salam}, F.~A., May 2014.
  {Relativistic and the first sectorial harmonics corrections in the critical
  inclination}. Astrophysics and Space Science 351, 113--117.

\bibitem[{{Sterne}(1958)}]{Sterne1958}
{Sterne}, T.~E., Jan. 1958. {The gravitational orbit of a satellite of an
  oblate planet}. The Astronomical Journal 63, 28.

\bibitem[{{Stone} and {Brodsky}(1988)}]{StoneBrodsky1988}
{Stone}, A.~D., {Brodsky}, R.~F., Aug. 1988. {Molniya orbits obtained by the
  two-burn method}. Journal of Guidance Control Dynamics 11, 372--375.

\bibitem[{{Ulybyshev}(2009)}]{Ulybyshev2009}
{Ulybyshev}, Y.~P., Aug. 2009. {Design of satellite constellations with
  continuous coverage on elliptic orbits of Molniya type}. Cosmic Research
  47~(4), 310--321.

\bibitem[{{Xu} and {Xu}(2013)}]{XuXu2013}
{Xu}, G., {Xu}, J., Feb. 2013. {On the singularity problem in orbital
  mechanics}. Monthly Notices of the Royal Astronomical Society 429~(1),
  1139--1148.

\bibitem[{Zhu et~al.(2014)Zhu, Zhao, and Zhang}]{ZhuZhaoZhang2014}
Zhu, T.-L., Zhao, C.-Y., Zhang, M.-J., 2014. Long term evolution of molniya
  orbit under the effect of earth's non-spherical gravitational perturbation.
  Advances in Space Research in press.
%\newline\urlprefix\url{http://www.sciencedirect.com/science/article/pii/S02731%
%17714002178}

\end{thebibliography}

\end{document}